\magnification\magstep1
\overfullrule = 0pt

\def\n{\noindent}
\def\qed{{\hfill{\vrule height7pt width7pt
depth0pt}\par\bigskip}} 
\def\vp{\varepsilon}
\def\pf{\medskip\n {\bf Proof.}~~}
\def\calt{{\cal T}}
\def\calc{{\cal C}}
\def\calk{{\cal K}}
\def\cale{{\cal E}}
\def\ms{\medskip}
\def\n{\noindent}
\def\ie{{\it i.e.\/}\ }
\def\cf{{\it cf.\/}\ }
\def\vp{\varepsilon}
\def\pf{\medskip\n {\bf Proof.}~~}
\def\calt{{\cal T}}
\def\calc{{\cal C}}
\def\calk{{\cal K}}
\def\cale{{\cal E}}
\def\cala{{\cal A}}
\def\calh{{\cal H}}
\def\ms{\medskip}
\def\olim{\limsup}
\def\ovl{\overline}
\def\EE{{\rm I}\!{\rm E}}
\def\CC{ \;Ê{}^{ {}_\vert }\!\!\!{\rm C}}
\def\RR{{\rm I}\!{\rm R}}
\def\ovl{\overline}
\def\EE{{\rm I}\!{\rm E}}
\def\CC{ \;Ê{}^{ {}_\vert }\!\!\!{\rm C}}

\centerline{\bf The similarity degree of an}
\centerline{\bf operator algebra, II}\bigskip
\centerline{ by Gilles Pisier\footnote*
{Supported in part by the NSF   and
 the Texas Advanced Research Program 010366-163.}}
\centerline{ Texas A\&M University}
\centerline{College Station, TX 77843, U. S. A.}
\centerline{and}
\centerline{Universit\'e Paris VI}
\centerline{Equipe d'Analyse, Case 186, 75252}
\centerline{ Paris Cedex 05, France}\bigskip

{\narrower\narrower\medskip\noindent {\bf Abstract.} For every integer $d\ge 1$, 
there is a unital closed subalgebra $A_d\subset B(H)$ with similarity degree 
equal precisely to $d$, in the sense of our
previous paper. This means that for  any unital
homomorphism $u\colon \ A_d\to B(H)$ we have
$\|u\|_{cb} \le  K\|u\|^d$ with $K>0$ independent
of $u$, and the exponent $d$ in this estimate 
cannot be improved. The proof that the degree
is larger than $d-1$  crucially uses an upper
bound for the norms of certain Gaussian random
matrices due to Haagerup and Thorbj\o rnsen. We
also include several complements to our previous 
publications on the same subject.\smallskip}
\bigskip
\centerline{\it Revised June 18, 1999.  
}\centerline{\it   To appear in
Math. Zeit.}
\vfill\eject

\baselineskip = 18pt

\n {\bf \S 0. Introduction}

This article is a continuation of our earlier papers [P1, P2]. We denote by 
$B(H)$ the algebra of all bounded operators on a Hilbert space $H$. Let $A$ be a 
unital operator algebra \ie a closed unital
subalgebra of $B(H)$. Assume that  every {\it
bounded\/} morphism (= unital homomorphism)
$u\colon \ A\to B(H)$ is  automatically
completely bounded (c.b.\ in short). Then (\cf  
[P1]) there is an  integer $d$ and a constant $K$
such that any such $u$ satisfies
$$\|u\|_{cb} \le K\|u\|^d.$$
The smallest $d$ for which this holds is called the similarity degree of $A$ and 
is denoted by $d(A)$. By convention, we set $d(A) = \infty$ if there is a 
bounded morphism $u\colon \ A\to B(H)$ which is not c.b. The main result proved 
in section~2 below is

\proclaim Theorem 0.1. For any $d\ge 1$, there is a (nonself-adjoint) unital 
operator algebra $A_d$ such that $d(A_d) = d$.

 Note that the existence of 
unital (nonself-adjoint) operator algebras $A$ with $d(A) = \infty$ is well 
known. In view of this, the preceding result is not too surprising. However its 
verification has proved to be much more difficult
than expected, although the  algebras $A_d$
themselves are rather canonical and easy to
define.

In [P1], we gave examples of $C^*$-algebras with degree equal to 1, 2 and 3 but 
we could not construct any examples (self-adjoint or not) with {\it finite\/} 
degree $>3$. The preceding result fills this gap in the nonself-adjoint case, 
but the case of $C^*$-algebras remains open. Note that a well known conjecture 
of Kadison [Ka] implies, modulo [P1], that there
is a universal bound for the  similarity degree
of $C^*$-algebras, but we are convinced that the
opposite is  what happens.

 The proof of Theorem~0.1 uses ``maximal operator
spaces'' in the  sense of [BP]. A typical
example is the space $\ell_1$ (or its
$m$-dimensional  version $\ell^m_1$) equipped
with the dual operator space structure to the 
commutative $C^*$-algebra $c_0$ (or
$\ell^m_\infty$). The next result, perhaps of
independent interest, is the main key tool which
we  use to prove Theorem~0.1.

\proclaim Theorem 0.2. Let $m\ge 1$ and $d\ge 1$ be integers. We denote simply 
$\ell^m_1$ for $\max(\ell^m_1)$. 
Let $E^d_m =  \ell^m_1\otimes_h\cdots \otimes_h 
\ell^m_1$ ($d$ times) and let
$$J_{m,d}\ \colon \ E^d_m\to
\max(\ell^m_2\otimes_2\cdots \otimes_2 
\ell^m_2)$$
be the identity map. Then
$$2^{-2(d-1)} m^{d-1\over 2} \le \|J_{m,d}\|_{cb} \le m^{d-1\over 2}.\leqno 
(0.1)$$
Moreover, the identity map $Id$ on $E^d_m$
satisfies
$$   2^{-2(d-1)} m^{d-1\over 2} \le \|Id :
E^d_m\to \max(E^d_m)\|_{cb}
\le m^{d-1\over 2}.\leqno  (0.2)$$

We can restate the preceding result in non-technical terms. In order to do this, 
we will now give a more explicit description of the constants estimated in (0.1) 
and (0.2) above.

Consider a family of scalars $\{\lambda_i\mid i= (i_1,\ldots, i_d) \in 
[1,\ldots, m]^d\}$ and let   $$P = \sum \lambda_i
X^1_{i_1}X^2_{i_2}\ldots  X^d_{i_d}$$ be the associated
``polynomial'' (homogeneous of degree $d$) in the
non-commutative variables 
$\{X^k_j\mid 1\le k \le d, 1\le j \le m\}$. We introduce the norm
$$|||P||| = \sup \left\|\sum \lambda_i x^1_{i_1} x^2_{i_2}\ldots 
x^d_{i_d}\right\| _{B(\ell_2)}$$
where the supremum runs over all families $\{x^k_j\mid 1\le k \le d, 1\le j \le 
m\}$ in the unit ball of $B(\ell_2)$ (actually,
the supremum is the same
if we restrict to the case when
$x^1_j=x^2_j=...=x^d_j$ for all $1\le j\le m$.).

We denote by $C(m,d)$ the smallest constant $C\ge 0$ with the
following property:\  for any family $\{x_i\mid i\in [1,\ldots,
m]^d\}$ in $B(\ell_2)$ such that
$$
\left\|\sum
\lambda_ix_i\right\|\le
\left|\left|\left|\sum
\lambda_i X^1_{i_1}  X^2_{i_2}\ldots
X^d_{i_d}\right|\right|\right| \leqno (0.3)\quad {\forall
\lambda_i
\in \CC  \quad (i\in [1,\ldots, m]^d)}$$   we 
can find operators\allowbreak $\{x^k_p\mid 1\le k
\le d, 1\le p\le m\}$ in the unit  ball of 
$B(\ell_2)$ such that
$$x_i = C \ x^1_{i_1} x^2_{i_2}\ldots x^2_{i_d}.\leqno
(0.4)\qquad \forall  i\in [1,\ldots, m]^d $$
(Note: Let $(e_i)$ be the natural basis of 
$E^d_m$. If we
fix the family
$(x_i)$ then the
 smallest constant $C\ge 0$ for which 
(0.4) holds is the $cb$-norm of the map 
from $E^d_m$ to $B(\ell_2)$ taking $e_i$ 
to $x_i$ and (0.3) means that this map is assumed of norm $\le 1$.) It is easy 
to check that
$$\left(\sum |\lambda_i|^2\right)^{1/2} \le \left|\left|\left| \sum \lambda_i 
X^1_{i_1} \ldots X^d_{i_d}\right|\right|\right|$$
so that for (0.3) to hold it suffices to have
$$\left\|\sum \lambda_ix_i\right\| \le \left(\sum |\lambda_i|^2\right)^{1/2} 
.\leqno (0.5)\quad {\forall
\lambda_i
\in \CC  \quad (i\in [1,\ldots, m]^d)}$$
  We denote by $C'(m,d)$ the 
smallest constant $C$ such that any family $(x_i)$ in $B(\ell_2)$ satisfying the 
estimate (0.5) admits a factorization of the form (0.4).

Then by known results (see [CS2]), we have $C'(m,d) = \|J_{m,d}\|_{cb}$ and 
$C(m,d) = \|Id\colon \ E^d_m\to
 \max(E^d_m)\|_{cb}$, whence the following 
reformulation of Theorem~0.2.

\proclaim Theorem 0.2 bis. We have for all $m,d\ge 1$
$$2^{-2(d-1)} m^{d-1\over 2} \le C'(m,d) \le C(m,d) \le m^{d-1\over 2}. \leqno 
(0.6)$$

Note that when the ``degree'' $d$ is fixed and $m\to \infty$ our estimates give 
the sharp order of magnitude (i.e.\ the exponent $(d-1)/2$ is sharp). The 
delicate point is the lower bound in (0.1), (0.2) or (0.6). Our proofs of this 
uses random matrices. In particular, we make crucial use of a remarkable upper 
bound for the norm of Gaussian random matrices with matrix coefficients due to 
Haagerup and Thorbj\o rnsen [HT] which we state below as Theorem~1.1. Indeed, 
the operators $x_i$ which achieve the lower bound in (0.6) are
actually obtained  by a random selection of matrices of
suitably large size which we can always  view as elements of
$B(\ell_2)$ by adding coefficients equal to zero. However,  the
precise form of our random matrices is rather complicated (see
(2.4) below).

As obvious idea which comes to mind is to let $\{x_i\mid i\in [1,\ldots, m]^d\}$ 
be an independent collection of $N\times N$ random matrices with independent 
entries, each one being Gaussian with mean zero and $L_2$-norm equal to 
$N^{-1/2}$. It is well known that, up to a numerical factor, when $N\to \infty$ 
(and $m,d$ remain fixed), these will satisfy (0.5) with large probability (say 
$>1/2$). We are thus reduced to estimate 
the best possible $C$ for which (0.4) 
holds for these. However, although this works when $d=2$, this
choice of
$(x_i)$ is  definitely not the right one when $d$ is larger.
Indeed, we show at the end of 
\S 2 that this way of choosing $x_i$ leads to weaker lower bounds in (0.6) when 
$d>2$ (and even that the resulting estimate is not sensitive
enough to ``distinguish" between the cases $d=2k$ and
$d=2k+1$!). This explains why we use a more complicated
definition of our random  choice of the matrices
$(x_i)$, in which the above mentioned result from [HT] is 
crucial to show that (0.5) still holds.

The algebras $A_d$ appearing in Theorem~0.1 can be described as follows. Let $E$ 
be an operator space and let $OA(E)$ be the universal unital operator algebra 
generated by $E$. This means that any complete contraction $v\colon \ E\to B(H)$ 
uniquely extends to a completely contractive morphism $\hat v\colon \ OA(E)\to 
B(H)$. We construct $OA(E)$ as a suitable completion of the tensor algebra of 
$E$. Let $I_d(E) \subset OA(E)$ be the ideal generated by $E\otimes\cdots 
\otimes E$ $ ((d+1)-\hbox{times})$. Then we
define $A_d(E) = OA(E)/I_d(E)$, and we let
$A_d=A_d(\ell_1)$. By [BRS], we know that these
quotients are (completely isometrically) operator
algebras. The  estimates in Theorem~0.2 will allow
us to prove:

\proclaim Theorem 0.3. Let $E$ be any infinite dimensional maximal operator 
space. Then  $$
 d(A_d(E)) =~d \qquad{\forall d \ge 1}.$$

\n {\bf Remark.}  If $E$ is not completely
isomorphic to a maximal operator space, then
$d(A_d(E)) =\infty$. Indeed, this follows from
[Pa3], where the case
$d=1$ of the preceding result is proved. More
precisely,
note that if $k<d$, then $ A_k(E)  $ is a quotient
of $ A_d(E)  $, so that $d(A_k(E)) \le d(A_d(E))$.
Taking $k=1$, we find that if the degree of $
A_d(E)$ is finite, then necessarily
$d(A_1(E))<\infty$ which implies (see [Pa3])
that $E$ is completely isomorphic to a maximal
operator space since any bounded map $v: E \to
B(H)$ must be c.b.

\n {\bf Remark.} By the non-commutative version
of von Neumann's inequality proved in [Bo], it is
easy to see that $OA(\ell_1)$ can be identified
(completely isometrically) with the unital closed
subalgebra of
$C^*(F_\infty)$ generated by the generators only
(and {\it not } their inverses); here
$F_\infty$ denotes the free group with countably
infinitely many generators. Moreover, 
$OA(E)$ coincides with the unital closed
subalgebra
generated by
$E$ in the
$C^*-$algebra of
$E$ in the sense of [Pes].

 After some
background in
\S 1, we prove these results at the end of  \S
2. Then in
\S 3 we  prove several complements. We return to
the general framework adopted in [P1] of  a
similarity setting, \ie an operator  space
generating an operator algebra.  In particular,
we will prove:

\proclaim Theorem 0.4. Let $A$ be a unital
algebra. Assume that, for some 
$\vp>0$, any morphism $u\colon \ A\to B(H)$ with $\|u\|\le 1+\vp$ is completely 
bounded. Then necessarily the same property holds for all $\vp>0$ and hence 
$d(A) < \infty$.

More generally, for any $c\ge 1$, let us denote by $\calc_c$ the class of all 
morphisms $u\colon \ A\to B(H)$ such that $\|u\| \le c$.

In addition, we will say that two morphisms $u_1\colon \ A\to B(H_1)$ and 
$u_2\colon \ A\to B(H_2)$ are similar if there exists an isomorphism $\xi\colon 
\ H_1\to H_2$ such that
$$u_2(a) = \xi u_1(a)\xi^{-1}.\leqno \forall~a\in A$$
Then we can state one more result to be proved in \S 3.

\proclaim Theorem 0.5. Let $A$ be a unital
operator algebra and let $1\le \theta  < c <
\infty$ be fixed. Then the following are
equivalent:
\item{(i)} Every morphism in $\calc_c$ is similar to a morphism in 
$\calc_\theta$.
\item{(ii)} Every bounded morphism $u\colon \
A\to B(H)$ is similar to one in 
$\calc_\theta$.
\medskip

\n {\bf Remark 0.6.} When (i) and (ii) above hold, the results of [P1] can be 
applied and yield that there are $\alpha>0$ and $K$ such that for any bounded 
morphism $u\colon \ A\to B(H)$ there is an invertible operator $\xi\colon \ H\to 
H$ with $\|\xi\|\ \|\xi^{-1}\| \le K\|u\|^\alpha$ such that $\|u_\xi\|_{cb}\le 
\theta$, where we have set $u_\xi(\cdot) = \xi^{-1}u(\cdot)\xi$. Moreover, the 
smallest such $\alpha$ is an integer. This is nothing but the similarity degree 
of a certain ``enveloping operator algebra'' which is denoted $\tilde A_\theta$ 
in [P1].

\n{\bf Acknowledgement.} I am very grateful to
Marius Junge for useful related information and
to C. Le Merdy for stimulating conversations. 
\vfill\eject

\n {\bf \S 1. Background}

 We recall  that an
``operator space'' is a closed subspace $E
\subset B(H)$ of the 
$C^*$-algebra of all bounded operators on a Hilbert space $H$. When $H = 
\ell_2$, we will denote by $\calk$ the subalgebra of all compact operators on 
$\ell_2$.  Let $E_1,E_2$ be operator spaces. We denote by $E_1 \otimes E_2$ 
their algebraic tensor product (as vector spaces). Assume $E_i \subset B(H_i)$ 
$(i=1,2)$. Then $E_1\otimes E_2$ can be identified with a linear subspace of 
$B(H_1 \otimes_2 H_2)$. The completion of $E_1\otimes E_2$ for the induced norm 
is called the minimal (= spatial) tensor product and is denoted $E_1 
\otimes_{\rm min} E_2$. Obviously the resulting embedding $E_1\otimes_{\rm min} 
E_2 \subset B(H_1\otimes_2 H_2)$ allows us to view $E_1 \otimes_{\rm min} E_2$ 
as an operator space. We will denote its norm by $\|~~~\|_{E_1\otimes_{\rm min} 
E_2}$, or simply by $\|~~~\|_{\rm min}$ when there is no ambiguity.

If $\dim (H) = n$, we identify $B(H)$ with the space $M_n$ of all $n\times n$ 
matrices with complex entries equipped with the usual operator norm.

Then, if $E$ is an operator space, $M_n
\otimes_{\rm min} E$ can be identified with  the
space $M_n(E)$ of all $n\times n$ matrices with
entries in $E$. In  particular, if $E = M_p$ for
some integer $p\ge 1$, $M_n \otimes_{\rm min}
M_p$  can be identified with $M_{np}$. Let $I_X$
denote the identity on a space $X$.  A linear
mapping $u\colon \ E_1\to E_2$ is called
completely bounded (in short  c.b.) if
$I_{\calk} \otimes u$ defines a bounded linear
map from $\calk \otimes  _{\rm min} E$, to $\calk
\otimes_{\rm min} E_2$, and we set
$$\|u\|_{cb} = \|I_{\calk} \otimes u\colon \
\calk \otimes_{\rm min}E_1 \to 
\calk \otimes_{\rm min} E_2\|.$$
For short we will often write $I$ instead of
$I_{\calk}$. We refer the reader  to [Pa1] and
[P3] for more information on c.b.\ maps and to
[BP] and [ER1--ER3]  for more on ``Operator
Space Theory''.

A mapping $u$ with $\|u\|_{cb}\le 1$ 
is called ``completely contractive'' or a
``complete contraction'', which we both abbreviate
by c.c. We will also use the  abbreviations o.s.\
and o.s.s.\ for ``operator space'' and ``operator
space  structure''.

By the term ``morphism'' we always mean a unital homomorphism between two unital 
algebras. 

We will need the notion of ``sum'' of
operator spaces, in the sense of  [P5]. This is
defined as follows. Let $\{E_\pi\mid \pi \in I\}$
be a finite  family of operator spaces indexed by
some (finite) index set $I$.

\n We assume that this family is ``compatible''
(this is the term used in  interpolation theory,
\cf  [BL]) \ie we assume given a specific
family 
$J_\pi\colon \ E_\pi\to X$ of continuous linear injective maps into a common 
Banach space $X$. Thus we may think of the spaces
$E_\pi$ as ``included'' in $X$. This allows us to 
``compare'' an element
$x$ in
$E_\pi$  with one
$x'$ in
$E_{\pi'}$ $(\pi\ne \pi')$ by considering their
images in $X$.  Thus the Banach space
$\sum\limits_{\pi\in I} E_\pi$ is defined as the
subspace  of $X$ formed of all elements $x$ in
$X$ of the form $x = \sum\limits_{\pi\in I} 
J_\pi(x_\pi)$ with $x_\pi \in E_\pi$ for all
$\pi$. Equipped with the norm 
$\|x\| = \inf \sum \|x_\pi\|$ (with the infimum running over all possible 
representations of $x$), this space becomes a
Banach space. The latter space  can be identified
with the quotient
$$\left(\oplus \sum_{\pi\in I} E_\pi\right)_1\bigg/ \Delta$$
where $\left(\oplus \sum\limits_{\pi \in I} E_\pi\right)_1$ denotes the 
$\ell_1$-direct sum of the family and where
$$\Delta = \left\{(x_\pi)_{\pi\in I}\ \Big| \ \sum\limits_{\pi\in I} x_\pi = 
0\right\}.$$
By classical results from operator space theory
 (based on Ruan's Theorem), the notions of
$\ell_1$-direct sum and the notion of quotient
have been  extended from the Banach to the
operator space category. Therefore the same is 
true of  course for the above ``sum''
$\sum\limits_{\pi\in I} E_\pi$. For details,
 see [P5, \S 2] or [P6, p.35]. Note  that we
only use the case when the index set $I$ is
finite and in that case the  space $\calk
\otimes_{\rm min} \left(\sum\limits_{\pi \in I}
E_\pi\right)$ can  be identified with {\it
equivalent\/} norms with the space
$\sum\limits_{\pi \in  I} \calk \otimes_{\rm min}
E_\pi$, but for $I$ infinite this does not
remain  true.

We will use repeatedly the notion of ``maximal'' operator space introduced
in [BP], and further studied in [Pa2]. Let us
recall its definition:\ let
$E$ be any normed space. Let $I$ be the class of all maps $u\colon \ B\to
B(H_u)$ with $\|u\|\le 1$ (and say $\dim H_u \le
  {\rm card}(E)$). We let
$J\colon \ E\to \bigoplus\limits_{u\in I} B(H_u)$ be the isometric
embedding defined by $J(x) = \bigoplus\limits_{u\in I} u(x)$. Then,
$\max(E)$ is defined as the operator space $J(E) \subset
B\left(\bigoplus\limits_{u\in I}H_u\right)$, and any operator space which
is of this form (up to complete isometry) is called ``maximal''.

The ``maximal'' operator spaces are characterized by the property that,
for any linear map $u\colon \ E\to B(H)$ we
 have $\|u\|_{cb} = \|u\|$. The
following slightly more explicit description
 of their operator space
structure from [Pa2] is often useful:\ for any $n$
and any
$x$ in $M_n(\max(E))$ we have $\|x\|<1$ iff, for
some integer $N$, there is a diagonal matrix $D$
in $M_N(E)$ and scalar matrices $\beta\in
M_{n,N}$ and $\gamma\in M_{N,n}$ such that 
$$x = \beta D\gamma\quad {\rm and}\quad 
\|\beta\|\,
\|D\|\, \|\gamma\|<1. $$
We refer the reader to
[Pa2] for more information on this.

By a ``complex Gaussian'' random variable, we mean a $\CC$-valued random 
variable $g$ with mean zero such that its real and imaginary parts are 
independent Gaussian variables {\it with equal variance}, so that the covariance 
matrix of $g$ viewed as $\RR^2$-valued is a
multiple of the identity.

\n We will make crucial use of the following
 remarkable result of Haagerup and 
Thorbj\o rnsen.

\proclaim Theorem 1.1. {\bf ([HT])}\ Let $r\ge 1$ and let $g_1,\ldots, g_r$ be a 
collection of independent random $N\times N$ matrices such that the entries 
$(g_k)_{ij}$ are independent (mean zero) Gaussian complex valued random 
variables with $\EE|(g_k)_{ij}|^2 = 1/N$, whenever $1\le k\le r$, $1\le i,j 
\le N$. Then for any $p$ and any $a_1,\ldots, a_r$ in $M_p$ we have almost 
surely
$$\limsup_{N\to \infty}  \left\| \sum^r_{j=1}
a_j\otimes g_j\right\|_{M_p  ( M_N)}
  \le \left\|\sum
a^*_ja_j\right\|^{1/2} + \left\|\sum 
a_ja^*_j\right\|^{1/2},$$ hence a fortiori $\le
2\left(\sum \|a_j\|^2\right)^{1/2}$.\ms

\n {\bf Remark.} For example, let us consider a family $\{g_{ij}\}$ as in 
Theorem~1.1 but indexed by $r=m^2$ this time so that $g_{ij}$ is indexed by a 
pair $i,j$ with $1\le i,j\le m$. We denote as usual by $(e_{ij})$ the canonical 
basis of $M_m$ and we introduce the random matrix (of size $mN\times mN$)
$$Y_{m,N} = \sum^m_{ij=1} e_{ij} \otimes g_{ij}.$$
Observe then that
$$\left\|\sum_{ij} e_{ij}e^*_{ij}\right\| =  \left\|\sum_{ij}
e_{ij}^*  e_{ij}\right\| = m.$$
Hence Theorem~1.1 implies in particular that
$$\limsup
\limits_{N\to \infty} \|Y_{m,N}\| \le 2m^{1/2}. 
\leqno  (1.1)$$

Note that actually, 
 as observed in [HT], inequalities such as (1.1) 
are      known to probabilists (\cf [Ge]) and can
be obtained by a much more direct proof (not
using  [HT]) but since the remaining part of our
argument depends crucially on [HT],  for brevity
we content ourselves with the preceding
derivation of (1.1) from  Theorem~1.1.

\proclaim Corollary 1.2. Fix an integer $d$. Let $(g^{(1)}_1,\ldots, g^{(1)}_r)$ 
be random matrices as in Theorem~1.1 but with size $N_1\times N_1$. Then let 
$(g^{(2)}_1,\ldots, g^{(2)}_r)$ be an analogous $r$-tuple but with size $N_2 
\times N_2$ and independent of the preceding collection, and so on until we 
reach $(g^{(d)}_1,\ldots, g^{(d)}_r)$. Now fix an integer $p$ and let 
$a_{j_1\ldots j_d}$ be a collection in $M_p$ indexed by $(j_1,\ldots, j_d)\in 
[1,\ldots, r]^d$. We have then almost surely
$$\eqalign{&\limsup_{N_1\to
\infty}\left(\limsup_{N_2\to\infty} \ldots \left( 
\olim_{N_d\to\infty} \left\|\sum a_{j_1j_2\ldots j_d} \otimes g^{(1)}_{j_1} 
\otimes g^{(2)}_{j_2} \otimes g^{(d)}_{j_d}\right\|\right)\right)\cr
&\quad \le 2^d\left(\sum \|a_{j_1j_2\ldots j_d}\|^2\right)^{1/2}}$$
where the summations run over all indices $(j_1,\ldots, j_d)$ with $1\le j_k \le 
r$.

\pf This is an immediate consequence of Theorem~1.1, by iterated 
applications.\qed\ms

We now wish to estimate the following random variable 
$$Z(r,p; N_1,\ldots, N_d)(\omega) = \sup \left\{\left\|\sum a_{j_1\ldots j_d} 
\otimes g^{(1)}_{j_1}(\omega)\otimes\cdots \otimes g^{(d)}_{j_d}(\omega) 
\right\| \right\}$$
where the sup runs over all families $\{a_{j_1\ldots j_p}\mid 1\le j_k \le r\}$ 
in $M_p$ such that $\sum\|a_{j_1\ldots j_p}\|^2\le 1$.

We will prove

\proclaim Corollary 1.3. For all $p$ and $r$ we have almost surely
$$\olim_{N_1\to \infty} (\olim_{N_2\to \infty} \ldots (\olim_{N_d\to \infty} 
Z(r, p; N_1,\ldots, N_d))) \le 2^d.$$

\pf Let $X$ be the finite dimensional Banach space of all families $a = 
(a_{j_1\ldots j_d})$ with $a_{j_1\ldots j_d} \in M_p$ equipped with the norm 
$\|a\| = \left(\sum \|a_{j_1\ldots
j_d}\|^2\right)^{1/2}$. Let $B$ denote  the
(compact) unit ball of $X$. Fix $0<\vp<1$ and let
$\Lambda$ be a finite 
$\vp$-net in $B$. By a classical elementary argument (\cf
 e.g.\ [P4, p. 49]) we  have
$$Z(r, p; N_1,\ldots, N_d) \le (1-\vp)^{-1} \sup_{a\in\Lambda} \left\|\sum 
a_{j_1\ldots j_d} \otimes g^1_{j_1}\otimes\cdots \otimes g^d_{j_d}\right\|.$$
Now, since $\Lambda$ is finite, Corollary~1.2 implies that we have almost surely
$$\olim_{N_1\to \infty} \left(\ldots \left(\olim_{N_d\to \infty} 
\sup_{a\in\Lambda} \left\|\sum a_{j_1\ldots j_d} \otimes g^1_{j_1} 
\otimes\cdots\otimes g^d_{j_d}\right\|\right)\right) \le 2^d.$$
Then  Corollary~1.3 follows immediately since
$\vp$ can be chosen  arbitrarily small.
\qed\ms

We will also need the following elementary fact.

\proclaim Lemma 1.4. With the same notation as above, let us denote for each $j 
= (j_1,\ldots, j_d)$ in $[1,\ldots, m]^d$
$$U_j = g^1_{j_1} \otimes \cdots\otimes g^d_{j_d}.$$
We have then almost surely
$$\liminf_{N_1\to \infty} \left(\ldots \left(\liminf_{N_d\to \infty} \left\|\sum 
U_j \otimes \ovl U_j \right\|\right)\ldots\right) \ge m^d.$$

\pf Let $N = N_1\times\cdots \times N_d$. It is well known that
$$\left\|\sum U_j\otimes \ovl U_j\right\| = \sup\left\{\left|\hbox{tr}\left( 
\sum U_j x U^*_j y\right)\right|\right\}$$
where the supremum runs over $x,y$ in the unit ball of the $N\times N$ 
Hilbert-Schmidt matrices. Taking $x$ and $y$ both equal to
$N^{-1/2}$-times the  identity we obtain
$$\left\|\sum U_j\otimes\ovl U_j\right\| \ge N^{-1}
\hbox{tr}\left(\sum  U_jU^*_j\right). 
\leqno  {(1.2)}$$
But
$$N^{-1} \hbox{tr}(U_jU^*_j) = \prod^d_{k=1}
[N^{-1}_k  \hbox{tr}(g^k_{j_k}{g^{k}_{j_k}}^*)]$$
and, for each $k$, by the   law of
large numbers (and the concentration of the
$\chi^2$-distribution around its mean) we know
that almost surely
$$\lim_{N_k\to \infty} {1\over N_k} \hbox{tr}
(g^k_j{g^{k}_j}^*) = 1.$$ Hence, if
$\min(N_1,\ldots, N_d)\to \infty$ we have for
each $j$
$$N^{-1}\hbox{tr}(U_j\otimes \ovl U_j)\to 1$$
whence
$$N^{-1}\hbox{tr}\left(\sum U_j\otimes \ovl
U_j\right)
\to m^d,$$ from which Lemma~1.4 follows
immediately by (1.2). \qed
\vfill\eject

\n {\bf \S 2. The examples}

The algebras $A_d$ are somewhat canonical. To describe them we start with a 
universal algebra $OA(E)$ which can be defined as
follows.

Let $E$ be an operator space. Let $\calt(E)$ be the tensor algebra of $E$, 
\ie $\calt(E) = \CC \oplus E \oplus E^{\otimes 2} \oplus \cdots$~. Any element of 
$\calt(E)$ can be written as a finite sum $x = x_0 + x_1+\cdots$ with $x_d \in 
E^{\otimes d}$ for all $d\ge 1$ and $x_0\in \CC$. We will denote by $P_j\colon\ 
\calt(E)\to E^{\otimes j}$ the
 mapping defined by $P_jx = x_j$. Any linear map 
$v\colon \ E\to B(H)$ admits a unique canonical extension $\hat v\colon \ 
\calt(E)\to B(H)$, as a morphism on the
unital algebra
$\calt(E)$. Let $I$ be the  collection of all
linear maps $v\colon \ E\to B(H_v)$ with
$\|v\|_{cb}\le 1$ (and with, say, $\rm
card(H_v)\le
\rm card(E)$),  and let
$J\colon
\
\calt(E)\to \bigoplus\limits_{v\in I} B(H_v)$ be
the mapping  defined by
$$J(x) = \bigoplus_{v\in I} \hat v(x).$$
Using this embedding, we equip $\calt(E)$ with a unital operator algebra 
structure, and we denote by $OA(E)$ the completion of the latter. We have 
clearly a canonical completely isometric embedding $E\subset OA(E)$. More 
generally, it is known (see [P1, Prop.~1.10] for details) that the closed 
subspace of $OA(E)$ generated by $E^{\otimes d}$ can be identified with the 
Haagerup tensor product $E\otimes_h \cdots \otimes_h E$ ($d$ times). Let us 
denote this subspace by $E_d$. For any $d\ge 1$, we denote by $I_d$ the closed 
ideal generated in $OA(E)$ by $E_{d+1}$. Equivalently $I_d$ can be described as 
the closed span of $\{E_m\mid m>d\}$. We can then form the quotient algebra
$$A_d(E) = OA(E)/I_d.$$
In the particular case $E = \ell_1$, we simply denote it by
$$A_d = A_d(\ell_1) = OA(\ell_1)/I_d.$$
We will need to describe a bit more the structure of the space $A_d(E)$, as 
follows.

\proclaim Proposition 2.1.
\item{(i)} $P_j\colon \ \calt(E)\to E^{\otimes j}$ defines a completely 
contractive projection on $OA(E)$.
\item{(ii)} Let $q\colon \ OA(E) \to A_d(E)$ be the canonical quotient map, and 
let $\widetilde E_j = q(E_j)$ $(0\le j\le~d)$. Then $\widetilde E_j$ is closed 
and
$$A_d(E) \simeq \widetilde E_0 \oplus \cdots \oplus \widetilde E_d.\leqno 
(2.1)$$
More precisely, any $x$ in $A_d(E)$ can be uniquely written as a sum $x = 
\sum^d_0 x_j$ with $x_j \in \widetilde E_j$, and the projection $x\to x_j$ is a 
complete contraction from $A_d(E)$ onto $\widetilde E_j$. 
\item{(iii)} Finally, for any $0\le j \le d$, the restriction $q_{|E_j}$ is a 
complete isometry from $E_j$ to $\widetilde E_j$.

\pf For any $z\in \CC$ with $|z|\le 1$, let $T(z)\colon \ OA(E)\to OA(E)$ be the 
completely contractive morphism associated to the linear map $v\colon \ E\to 
OA(E)$ equal to $z$ times the embedding of $E$ into $OA(E)$. For any $y$ in 
$E^{\otimes d}$, we have $T(z)y = z^dy$. In particular $T(z) (I_d) \subset I_d$, 
so that $T(z)$ canonically induces a completely contractive morphism 
$\theta(z)\colon \ A_d(E) \to A_d(E)$. Consider $y$ in $OA(E)$ of the form $y = 
\sum\limits_{j\ge 0} y_j$ (finite sum) with $y_j \in E_j$. We have $T(z)y = 
\sum\limits_{j\ge 0} z^jy_j$ hence $y_j = \int \bar zT(z)y
dm(z)$, which  implies $\|y_j\| \le \|y\|$.
Thus, we have a contractive projection $y\to
y_j$  from $OA(E)$ onto $E_j$. This averaging
argument actually shows that it is a  complete
contraction, whence (i). 

\n Applying the quotient
map $q$, we obtain  similarly
$$\theta(z) q(y) = \sum_{j\ge 0} z^j q(y_j)$$
whence again
$$\|q(y_j)\| \le \|q(y)\| = \left\|\sum q(y_j)\right\|$$
which  shows that the mapping $q(y)\to q(y_j)$ defines a completely contractive 
projection from $A_d(E)$ onto $\widetilde E_j = q(E_j)$. Also note that 
$\widetilde E_j$ is closed. Of course we have $q(E_j) = 0$ $\forall~j >d$, 
whence the decomposition (2.1), completing the proof of (ii).

\n Note that $q$ restricted to $E_0+\ldots+E_d$
is clearly injective. To show the last point
(iii), consider $y$ in $E_j$ ($0\le j\le d$) with
$\|q(y)\|<1$. Then there is $y'$ in $I_d$ with
$\|y+y'\|<1$.  But since
$P_j(y+y')=y$ we have $\|y\|<1$. Thus by
homogeneity we have
$\|y\| \le \|q(y)\|,$
which shows that $q_{|E_j}$ is an isometry. The proof that it is a complete 
isometry is similar and left to the reader.\qed

\proclaim Lemma 2.2. Let $E$ be a maximal operator space. Then any bounded 
morphism $u\colon \ A_d(E)\to B(H)$ is c.b.\ and satisfies
$$\|u\|_{cb} \le \sum\nolimits^d_0 \|u\|^j \le (d+1) \|u\|^d,$$
in particular $d(A_d(E))\le d$.

\pf Let $\calk_0 = \bigcup\limits_n M_n$. 
Consider an element $y$ in $\calk_0 
\otimes A_d(E)$ with $\|y\|<1$. We can write
$y=q(x) = q(x_0 +\cdots+ x_d)$ with 
$\|x\|<1$ and $x_j \in
\calk_0 
\otimes E_j$ where $q\colon \ \calk_0\otimes OA(E) \to \calk_0 \otimes A_d(E)$ 
is the quotient map. Therefore if we let $v =
uq_{|E}$ (we identify $E$ with 
$E_1$) we have
$$(I\otimes u)(y) = (I\otimes \hat v) [x_0
+\cdots+ x_d] =
\sum^d_0 (I\otimes v^{\otimes  j})(x_j)$$
hence
$$\eqalign{\|(I\otimes u)(y)\| &\le
\sum\nolimits^d_0 \|v^{\otimes j} \colon \
E_j\to  B(H)\|_{cb} \|x_j\|\cr &\le
\sum\nolimits^d_0 \|v\|^j_{cb}\|x_j\|}$$ and by
(i) in Proposition~2.1 we have $\|x_j\|\le
\|x\|<1$, whence the conclusion  since our
assumption that $E$ is ``maximal'' ensures that
$\|v\|_{cb}\le \|v\|\le 
\|u\|$.\qed\ms

Let $E$ be an operator space. Let $E_m$ be as above $(m\ge 0)$ with the 
convention that $E_0 = \CC$. 
For any partition $\pi = (m_1,\ldots, m_K)$ with 
$m_i\ge 0$ such that $\sum m_i = d$ we have a
natural product map
$$E_{m_1} \otimes \cdots\otimes E_{m_K}\to E_d$$
which is a complete contraction from
$E_{m_1}\otimes_h\cdots \otimes_h E_{m_K}$ to 
$E_d$. A fortiori, we have a complete contraction from
$$  \max(E_{m_1}) \otimes_h\cdots\otimes_h
\max(E_{m_K}) \hbox{ to }  E_d.$$
Let us denote whenever $\pi = (m_1,\ldots, m_K)$
$$E(\pi) =
 \max(E_{m_1}) \otimes_h\cdots\otimes_h
\max(E_{m_K}) .$$

If $m_i>0$ for each $i$, it can be shown 
that the latter map is always injective (we skip the details, see e.g.
[BP, Th. 3.11] which implies that the Haagerup
tensor
product of two {injective} maps is injective). This is
obvious
 if
$E$ is finite dimensional, or say if $E = \ell_1$.

\n Using these natural continuous injections
$E(\pi)\subset E_d$, we define
$$ \cale_k=  \sum_\pi E(\pi) $$
where the sum (this is {\it a  sum of operator
spaces} as defined above, in \S 1) is
relative to all partitions
$\pi = (m_1,\dots, m_K)$ of $d$ into 
$K$ disjoint nonempty blocks with $K\le k$ and
with
$\sum m_i = d$. For example, $\cale_1=\max(E_d)$,
$\cale_2=\max(E_d) +
\sum_{p=1}^{d-1}\max(E_p)\otimes_h \max(E_{d-p})
$, and so on.

\n Thus, by what precedes,  we have a complete
contraction
$$\Phi_k\colon \ \cale_k \to E_d.$$
\proclaim Lemma 2.3. If $d(A_d(E)) \le k$, then
necessarily $\Phi_k$ is a  complete isomorphism,
\ie $(\Phi_k)^{-1}$ is a c.b.\ map from $E_d$
to 
$\cale_k=\sum\limits_\pi E(\pi)$.

\pf It is proved in [P1, Th. 4.2] that if
$d(A_d(E))\le k$, then the product map  defines a
complete surjection from
$$\max(A_d(E)) \otimes_h\cdots \otimes_h \max(A_d(E))\quad (k \hbox{ times})$$
onto $A_d(E)$.

\n More precisely, let
$$X_k = \max(A_d(E)) \otimes_h\cdots \otimes_h
\max(A_d(E))\quad (k \hbox{  times})$$
and let us denote by $p_k\colon \ X_k\to A_d(E)$ this product map. Then there is 
a constant $C$ such that for any $x$ in $\calk \otimes_{\rm min} A_d(E)$ there 
is an element $y$ in $\calk \otimes_{\rm min} X_k$
with $\|y\| \le C\|x\|$ such that
$(I\otimes p_k) (y) = x$. By
Proposition~2.1, we have seen  that 
$$A_d(E) \simeq E_0 \oplus \cdots\oplus E_d$$
 hence
$$\max(A_d(E)) \simeq \max (E_0) \oplus \cdots\oplus \max(E_d)$$
and therefore
$$X_k \simeq \bigoplus[\max(E_{m_1}) \otimes_h \cdots \otimes_h \max(E_{m_k})]$$
where the direct sum runs over all 
$\pi = (m_1,\ldots, m_k)$ with $0\le m_i \le 
d$.

\n Moreover the projections onto the coordinates
of this direct sum are complete  contractions.
 Therefore $y$ can be written as
$$y = \sum_\pi y_\pi \quad \hbox{with}\quad y_\pi \in \calk \otimes E(\pi).$$
If we apply $T(z) \otimes\cdots\otimes T(z)$ ($k$ times) to this equality we 
obtain
$$(I\otimes T(z)^{\otimes k})y = \sum_\pi
z^{|\pi|} y_\pi$$ where $|\pi| = m_1+\cdots +m_k$.

\n Thus if we actually apply all this to an
element
$x$ in $\calk \otimes 
\widetilde E_d$, we have $x = (I\otimes p_k) (y)$ hence
$$\eqalign{  z^d x &= (I\otimes \theta(z))(x)\cr
&= (I \otimes \theta(z)p_k)(y)\cr
&= (I\otimes p_k T(z)^{\otimes k}) (y)\cr
&= \sum_\pi z^{|\pi|} (1\otimes p_k)(y_\pi).}$$
Therefore we must have
$${x = \sum_{\pi :|\pi|=d} (1\otimes
p_k)(y_\pi).}$$ Moreover since $y \to y_\pi$ is a
(completely) contractive projection, we have
$$ \sum _{\pi\colon \ |\pi|=d} \| y_\pi \|
\le C'\|y\| \le C'C\|x\|$$ where $C' =
\hbox{card}\{\pi\mid |\pi|=d\}$.
This is what we want, except that we have allowed
$m_i=0$ in
$\pi = (m_1,\ldots, m_k)$. Here is how this point
can be fixed: by deleting 
the zero blocks (but maintaining the order) each 
such $\pi$ defines a partition
$\hat\pi = (m'_1,\ldots, m'_K)$ with $m'_i
>0$ and $K\le k$. Then we can write
$$ \sum _{\pi\colon \ |\pi|=d} y_\pi =\sum_\omega
\sum_{\pi:\hat\pi=\omega} y_\pi $$
 where the sum $\sum_\omega$  is restricted to
partitions of
$d$ of the form
$\omega=(m'_1,\ldots, m'_K)$ with $m'_i
>0$ as appearing in the definition of
$\cale_k$.

\n Modulo   repeated identifications of
$\CC\otimes E$ with
$E$, we have 
$$x=(I\otimes \Phi_k)(\sum_\omega
y_\omega')\quad
{\rm where}\quad
y_\omega'=\sum_{\pi:\hat\pi=\omega} y_\pi
\in
\calk\otimes E(\omega).
$$ 
Thus we obtain that
$$\|[I\otimes \Phi^{-1}_k](x)\| \le C'C\|x\|$$
or equivalently $\|\Phi^{-1}_k\|_{cb} \le C'C$.\qed

\proclaim Lemma 2.4. When $k = d-1$ and $E$ is a maximal operator space, the 
space $\cale_k$ actually reduces (completely
isometrically) to
$$\sum_{\scriptstyle p+q+r=d\atop
 \scriptstyle q>1} E_p \otimes_h \max(E_q) 
\otimes _h E_r,\leqno (2.3)$$
with the convention that $E_p \otimes_h \max(E_q) 
\otimes _h E_r$ should be replaced
  by $E_p \otimes_h \max(E_q) 
$ if $r=0$ and $p>0$, by $\max(E_q) 
\otimes _h E_r$ if $p=0$ and $r>0$, and finally by
$\max(E_q) $ if $p=r=0$.

\pf Indeed, let $\sigma = (m'_i)$ with $m'_1 = \cdots= m'_p=1$,
$m'_{p+1}=q>1$ and 
$m'_{p+1+i} = 1$ for $i=1,\ldots, r$. Then
$E(\sigma) = E_p\otimes_h \max(E_q) 
\otimes_h E_r$ where of course $|\sigma| = \sum
m'_i = d$ and $\sigma$ has at most 
$d-1$ disjoint blocks since $q>1$. Thus when $k=d-1$, the sum appearing in 
(2.3) is canonically included in $\cale_k$. Now
take any $\pi = (m_i)_{i\le K}$  where $K<d$.
Since we have less than $d$ blocks we have
$m_i>1$ for some $i$.  Clearly if $p = m_1
+\cdots+ m_{i-1}$, if  $q=m_i$ and if $r =
m_{i+1} +\cdots+ m_{K}$  we have then
$$\max(E_{m_1}) \otimes_h \cdots \otimes_h \max(E_{m_K}) \subset E_p \otimes_h 
\max(E_q) \otimes_h E_r$$
whence the conclusion of Lemma~2.4.\qed

\proclaim Lemma 2.5. For any $m\ge 1$, we consider the operator space 
$\max(\ell^m_1)$. We will denote
$$E^d_m = \max(\ell^m_1)\otimes_h \cdots \otimes_h \max(\ell^m_1) \quad (d 
\hbox{ 
times}),$$
and also
$$\cale^{d-1}_m = \sum_{\scriptstyle p+q+r=d \atop \scriptstyle q>1} E^p_m 
\otimes_h \max(E^q_m) \otimes_h E^r_m$$
with the convention $E^0_m = \CC$. Then there is a constant $\delta>0$ 
(independent of $m$) such that for all $m\ge 1$ we have
$$\|i\colon \ E^d_m \to \cale^{d-1}_m\|_{cb}\ge \delta\sqrt m$$
where $i$ denotes the identity map.

  For the proof we need separate estimates as 
follows.

\proclaim Sublemma 2.6. Let $X$ be an arbitrary operator space and let 
$\otimes^{\wedge}$ denote the operator space projective tensor product. Then the 
identity map from $\ell^m_1 \otimes_h X$ (resp.\ $X \otimes_h \ell^m_1)$ into 
$\ell^m_1 \otimes^{\wedge} X$ (resp. $X \otimes^\wedge \ell^m_1$) has c.b.\ norm 
$\le \sqrt m$.

\pf The identity of $\ell^m_1$ factorizes through the ``row'' space $R_m$ as 
$\ell^m_1 {\buildrel a\over \longrightarrow} R_m {\buildrel a^{-1}\over {\hbox 
to 20pt{\rightarrowfill}}} \ell^m_1$ with $\|a\|_{cb}\le 1$, $\|a^{-1}\|_{cb}\le 
\sqrt m$. On the other hand it is well known (\cf
 [ER2, Th. 4.3] and [B, Prop. 2.3 (ii)] ) that
$R_m 
\otimes_h X$ can be identified with $R_m\otimes^\wedge X$ completely 
isometrically. Thus we obtain a factorization
$$\ell^m_1 \otimes_h X {\buildrel a\otimes I\over {\hbox to 
30pt{\rightarrowfill}}} R_m \otimes_h X = R_m \otimes^\wedge X {\buildrel 
a^{-1}\otimes I\over {\hbox to 30pt{\rightarrowfill}}} \ell^m_1 \otimes^\wedge 
X$$
from which the announced result follows immediately. The proof of the transposed 
statement is analogous (with $C_m$ instead of $R_m$).\qed

\n {\bf Remark.} More generally, let $Y$ be an $m$-dimensional operator space 
and let $a\colon \ R_m\to Y$ be a complete isomorphism. Then the preceding 
argument shows that
$$\|Y\otimes_h X\to Y \otimes^\wedge X\|_{cb} \le \|a\|_{cb} \|a^{-1}\|_{cb} 
\/.$$

\proclaim Sublemma 2.7. The identity map $i$ satisfies
$$\|i\colon \ \cale^{d-1}_m\to \max(E^d_m)\|_{cb} \le m^{d-2\over 2}.$$

\pf By the canonical property of the ``sum'', it suffices to show that
if $p+q+r=d$  with $q>1$ we have
$$\|i\colon \ E^p_m \otimes_h \max(E^q_m) \otimes_h E^r_m \to \max(E^d_m) \| 
_{cb} \le m^{d-2\over 2}.$$
Let $X = \max(E^q_m)$.
By iterated applications of Sublemma~2.6 (using
the associativity of $\otimes_h$  and
$\otimes^\wedge$, \cf  [BP, ER1-4]) we have
$$\|i\colon \ E^p_m \otimes_h X \otimes_h E^r_m\to \ell^m_1 \otimes^\wedge 
\cdots \otimes^\wedge \ell^m_1 \otimes^\wedge X \otimes^\wedge \ell^m_1 
\otimes^\wedge \cdots \otimes \ell^m_1\|_{cb} \le m^{p+r\over 2}.$$
Then we note that by the projectivity of $\otimes^\wedge$ the space $\ell^m_1 
\otimes^\wedge \cdots\otimes^\wedge \ell^m_1 \otimes^\wedge X \otimes^\wedge 
\ell^m_1 \cdots \otimes^\wedge \ell^m_1$ is clearly a maximal operator space  
(completely) contractively included in $E^d_m$.
(Indeed, by [BP, p. 289]  the
class of maximal operator spaces is stable under
the o.s. projective tensor product.) Since it is
maximal, it is  completely contractively included
in
$\max(E^d_m)$. Thus we obtain the announced 
result by noting simply that $q>1$ ensures
$p+r\le d-2$.\qed 
The next estimate is the key point of this paper.

\proclaim Sublemma 2.8. The identity map $i$ satisfies 
$$\|i\colon \ E^d_m\to \max(E^d_m)\|_{cb} \ge \delta m^{d-1\over 2}$$
where $\delta>0$ is a constant independent of $m$. (We will obtain $\delta = 
2^{-2(d-1)}$.)

\pf Fix $m\ge 1$. We denote $[m] = [1,2,\ldots, m]$. For any $i = (i_1,\ldots, 
i_d)$ in $[m]^d$, we denote by
$$e_i = e_{i_1} \otimes e_{i_2} \otimes \cdots \otimes e_{i_d}$$
the canonical basis vectors in the space $E^d_m$. Now we fix a number $0 < \vp 
< 1$ throughout the proof.

We will show that we can find matrices $\{U_i\mid i= (i_1,\ldots, i_d) \in 
[m]^d\}$ of size $N\times N$ (for a suitable $N =
N(m,d,\vp)$ such that:
\item{(i)} $\left\|\sum U_i \otimes \ovl U_i\right\| \ge (1-\vp) m^d$.
\item{(ii)} $\forall~\lambda_i \in \CC$
\ $(i\in[m]^d)$
 \qquad{$\left\|\sum \lambda_i U_i\right\| \le
(1+\vp) 2^{d-1} \left(\sum 
|\lambda_i|^2\right)^{1/2}$.}
\item{(iii)} The linear map $v\colon \ E^d_m \to M_N$ defined by $v(e_i) = U_i$ 
satisfies 
$\|v\|_{cb} \ge (1+\vp)^{-1} 2^{-d+1} m^{d-1\over
2}$.

\n  Let us first   check that (i), (ii) and
(iii) imply Sublemma~2.8. To do this, first
observe  that (ii) implies $\|v\| \le
(1+\vp)2^{d-1}$. This follows form the fact that
if 
$C_m$ denotes the column Hilbert space  in
dimension
$m$, then we have trivially a  complete
contraction $\ell^m_1 \to C_m$, hence a complete
contraction $E^d_m\to  C_m \otimes_h
\cdots\otimes_h C_m$ ($d$ times) but $C_m
\otimes_h \cdots 
\otimes_h C_m$ is (completely) isometric to the
$m^d$-dimensional column  Hilbert space (\cf 
[ER2, p. 272]). Hence (ii) implies $\|v\| \le
(1+\vp)2^{d-1}$.  Therefore, we have
$$\|i\colon \ E^d_m\to \max(E^d_m)\|_{cb} \ge {\|v\|_{cb}\over \|v\|} \ge 
(1+\vp)^{-2} 2^{-2(d-1)} m^{d-1\over 2},$$
and Sublemma~2.8 follows. Thus it suffices to
produce $(U_i)$ satisfying (i),  (ii), (iii).

The matrices $U_i$ will be defined as products of the following form
$$U_{i_1i_2\ldots i_d} = U^1_{i_1i_2} U^2_{i_2i_3}\ldots U^{d-1}_{i_{d-1}i_d},$$
where $U^k_{ij}$ are matrices in $M_N$ for all $1\le i,j\le m$. We will make 
sure that these matrices satisfy
\item{(iv)} For each $k=1,2,\ldots, d-1$,
$$\left\|\sum_{ij\le m} e_{ij}\otimes U^k_{ij}\right\|_{M_m(M_N)} \le 
2(1+\vp)m^{1/2}.$$

Let us now verify that (i) and  (iv)
together imply   (iii). Let us denote
by 
$\xi_i = \xi_{i_1} \otimes\cdots \otimes \xi_{i_d}$
 the canonical basis vectors 
of 
$(E^d_m)^* = \ell^m_\infty
\otimes_h\cdots\otimes_h \ell^m_\infty$. Let
$U^k  = \sum\limits^m_{p,q=1} e_{pq} \otimes
U^k_{pq}$ and let $\hat e_{pq} = e_{pq} 
\otimes I$. Then we have
$$e_{11} \otimes U_{i_1i_2\ldots i_d} = \hat e_{1i_1} U^1 \hat e_{i_2i_2} 
U^2\ldots \hat e_{i_{d-1}i_{d-1}} U^{d-1}\hat e_{i_d1}.$$
Let $V\colon \ (E^d_m)^*\to M_N$ be the linear map associated to the tensor 
$\sum U_i \otimes e_i \in M_N\otimes E^d_m$ so that $V(\xi_i) = U_i$. Then we 
have $\left\|\sum U_i \otimes e_i\right\|_{M_N(E^d_m)} = \|V\|_{cb}$, but the 
last identity implies that
$$V(\xi_{i_1} \otimes \cdots \otimes \xi_{i_d}) = v_1(\xi_{i_1}) v_2(\xi_{i_2}) 
\ldots v_d(\xi_{i_d})$$
where the linear maps $v_k \colon \ \ell^m_\infty\to M_m(M_N)$ are defined as 
follows, for all $1\le j\le m$
$ v_1(\xi_j)  = \hat e_{1j} U^1$, then 
$v_2(\xi_j) = \hat e_{jj}U^2$,
$ \dots$,
$v_{d-1}(\xi_j)  = \hat  e_{jj}U^{d-1}$
and finally 
$v_d(\xi_j)  = \hat e_{j1}.$
Then we have obviously for all $2\le k \le d-1$
$$\|v_k\|_{cb} \le \|U^k\|$$
hence using (iv) $\|v_k\|_{cb} \le 2(1+\vp)m^{1/2}$. Similarly, we have (by a 
well known fact)
$$\|v_1\|_{cb} \le m^{1/2}\|U^1\| \le 2(1+\vp)m$$
and
$$\|v_d\|_{cb} \le m^{1/2}.$$
Thus we conclude by the classical results on multilinear c.b.\ maps
([CS1-2, PS])  that we have 
$$\|V\|_{cb} \le \prod^d_1 \|v_k\|_{cb} \le (2(1+\vp))^{d-1} m^{d+1\over 2}.$$
Equivalently, this means that (\cf [BP, ER1])
$$\left\|\sum U_i\otimes e_i \right\|_{M_N(E^d_m)}
= \|V\|_{cb} \le  (2(1+\vp))^{d-1} m^{d+1\over
2}.$$ The same proof actually shows that\
$$\left\|\sum \ovl U_i 
\otimes e_i\right\|_{M_N(E^d_m)} \le
(2(1+\vp))^{d-1}  m^{d+1\over 2}.$$
But now by (i) and by the definition of the c.b.\ norm
(since $v(e_i)=U_i$) we have
$$\eqalign{(1-\vp)m^d &\le 
\left\|\sum U_i\otimes \ovl
U_i\right\|=\left\|\sum \ovl U_i\otimes 
U_i\right\|\cr &\le \|v\|_{cb(E^d_m,M_N)}
\left\|\sum \ovl U_i \otimes 
e_i\right\|_{M_N(E^d_m)}}$$ whence we find
$$\|v\|_{cb(E^d_m, M_N)} \ge (1-\vp) (2(1+\vp))^{-(d-1)} m^{d-1\over 2},$$
which concludes the proof that (i) and (iv)
imply  (iii).

\n  We now come to 
the main point:\ the construction of matrices $U_i$ satisfying (i), (ii) and 
(iv). 
Actually the product appearing in the definition of $U_i$ will be a tensor 
product, \ie we will have
$$U_{i_1i_2\ldots i_d} = Y^1_{i_1i_2} \otimes Y^2_{i_2i_3} \otimes\cdots \otimes 
Y^{d-1}_{i_{d-1}i_d}\leqno(2.4)$$
where $Y^k_{ij}$ are matrices of sizes $N_k\times N_k$.

In other words, we will set
$$U^k_{pq} = 1\otimes \cdots \otimes 1\otimes Y^k_{pq} \otimes 1 \otimes \cdots 
\otimes 1$$
and $N = N_1\times \cdots \times N_{d-1}$. The matrices $Y^k_{ij}$ will be 
chosen at random independently according to a Gaussian distribution. More 
precisely, the family $\{Y^k_{ij}\mid 1\le i,j\le m, 1\le k \le d-1\}$ is taken 
to be an independent collection of random variables, and for each $i,j$ and 
$k$ $(1\le i,j\le m, 1\le k\le d-1)$, $Y^k_{ij}$
is a random $N_k\times N_k$  matrix, the entries
of which $(Y^k_{ij})_{pq}$ $(1\le p,q\le N_k)$
are  independent complex valued Gaussian variables
with mean zero and such that 
$\EE|(Y^k_{ij})_{pq}|^2 = (N_k)^{-1}$, exactly as in Theorem~1.1 above. Then, 
for clarity of notation we introduce the random variables $Z_2$ and $Z_4$ 
defined by
$$\eqalignno{Z_2 &= \sup_{\sum |\lambda_i|^2 \le 1} \left\|\sum \lambda_i 
Y^1_{i_1i_2} \otimes\cdots \otimes Y^{d-1}_{i_{d-1}i_d}\right\|_{M_N}\cr
\noalign{\hbox{and}}
Z_4 &= m^{-1/2} \sup_{1\le k\le d-1} \left\|\sum^m_{p,q=1} e_{pq} \otimes 
Y^k_{pq}\right\|_{M_m(M_{N_k})}.}$$
Then by Corollary~1.3 we know that
$$\olim_{N_1\to \infty} \ldots \olim_{N_{d-1}\to \infty} Z_2 {\buildrel {\rm 
a.s.}\over \le} 2^{d-1}$$
and by Theorem~1.1 (and the remark following it) we have
$$\olim_{N_1\to \infty} \ldots \olim_{N_{d-1}\to \infty} Z_4 {\buildrel {\rm 
a.s.}\over \le} 2.$$
Since we also have almost surely
$$\liminf_{N_1\to \infty} \ldots \liminf_{N_d\to \infty} \left\|\sum U_j \otimes 
\ovl U_j\right\|\ge m^d,$$
it is now clear that the event corresponding to (i), (ii) and (iv) occurs
 with positive 
  probability (actually  
    $=1$) if $N_1, N_2,...,N_d$ are suitably
large, thus establishing the existence of
matrices satisfying  (i)--(iv) which completes
the proof of Sublemma~2.8.\qed \ms

\n {\bf Remark.} It is also possible to produce unitary matrices $(U_i)$ 
satisfying essentially the properties (ii), (iii), (iv) but with an additional 
numerical factor in front of the constants involved. (Hint:\ Use the 
concentration of measure phenomenon (see e.g.
[P4 , p. 44]) and a comparison principle such as 
the one appearing in [MP, p. 84].) The inequality
(i) is then automatically true  (and actually
becomes an equality).

\n {\bf Proof of Lemma 2.5.} We have obviously
$$\|i\colon \ E^d_m\to \max(E^d_m)\|_{cb} \le
\|i\colon \ E^d_m \to 
\cale^{d-1}_m\| _{cb}\cdot \|i\colon \ \cale^{d-1}_m \to \max (E^d_m)\|_{cb}.$$
Hence the result immediately follows from the preceding two sublemmas.\qed\ms

\n {\bf Proof of Theorem 0.1.} We let $E = \ell_1$ and $A_d = A_d(\ell_1)$. By 
Lemma~2.2, it suffices to prove that $d(A_d)>d-1$. Assume to the 
contrary that 
$d(A_d) \le  d-1$.  Then, by 
Lemma~2.3, $(\Phi_{d-1})^{-1}$ is a c.b.\ map from $E_d = \ell_1\otimes_h \cdots 
\otimes_h \ell_1$ ($d$-times) to $\cale_{d-1}$. By
Lemma~2.4, this means that we  have a c.b.\
mapping
$$E_d\to \sum_{\scriptstyle p+q+r=d\atop \scriptstyle q>1} E_p \otimes_h 
\max(E_q) \otimes_h E_r$$
which reduces to the identity on $E\otimes\cdots \otimes E$ ($d$ times). But 
this clearly contradicts Lemma~2.5.\qed

\n {\bf Proof of Theorem 0.2.}
By iterated applications of Sublemma~2.6, we find 
$$\|Id\colon \ E^d_m\to \ell^m_1 \otimes^\wedge\cdots \otimes^\wedge 
\ell^m_1\|_{cb}\le m^{d-1\over 2},$$
and (as already noted for Sublemma~2.7) we know that $\ell^m_1 \otimes^\wedge 
\cdots \otimes^\wedge \ell^m_1$ ($d$-times) is a maximal operator space (which 
can be identified with $\ell^{m^d}_1$), completely contractively included in 
$E^d_m$.
Therefore, we have a fortiori
$$\|Id\colon \ E^d_m\to \max(E^d_m)\|_{cb} \le m^{d-1\over 2},$$
and (as already observed in the proof of Sublemma~2.8) we have $\|J_{m,d}\|\le 
1$ whence
$$\|Id\colon \ \max(E^d_m)\to \max(\ell^m_2 \otimes_2\cdots \otimes_2 
\ell^m_2)\|_{cb}\le 1$$
hence the last estimate yields
$$\|J_{m,d}\|_{cb}\le m^{d-1\over 2}.$$
Thus we have proved the right side of both (0.1) and (0.2). The left side of 
(0.2) is but Sublemma~2.8. Now a close look at
the proof of Sublemma~2.8 shows  that the mapping
$v$ appearing there actually satisfies
$$\|v\colon \ \max(\ell^m_2 \otimes_2\cdots
\otimes_2
\ell^m_2)\to M_N\|_{cb}=
\|v\colon \ \ell^m_2 \otimes_2\cdots \otimes_2
\ell^m_2\to M_N\|\le  (1+\vp)2^{d-1}$$
and
$$\|v J_{m,d}\|_{cb}\ge (1+\vp)^{-1} 2^{-d+1} m^{d-1\over 2}.$$
Thus we must have $\|J_{m,d}\|_{cb} \ge (1+\vp)^{-2} 2^{-2(d-1)} m^{d-1\over 2}$ 
from which the left side of (0.1) follows.\qed

\n {\bf Proof of Theorem 0.3.} 
Let $E$ be a maximal  operator space
such that for each $m\ge 1$ the natural inclusion
$\ell^m_1\to 
\ell^m_2$ admits a factorization $\ell^m_1 {\buildrel a_m \over {\hbox to 
20pt{\rightarrowfill}}} E {\buildrel b_m \over {\hbox to 
20pt{\rightarrowfill}}} \ell^m_2$ with 
$$C=\sup\limits_m \|a_m\colon \ \max(\ell^m_1) \to
E\|_{cb} \| b_m\colon \  E\to 
\max(\ell^m_2)\|_{cb} < \infty.$$
We then claim that $d(A_d(E)) = d$.
 
\n Assume to the contrary that $d(A_d(E)) < d$.
Then as explained above, we must have a c.b.\
map  from
$E_d$ to
$\sum\limits_{\scriptstyle p+q+r=d\atop
\scriptstyle q\ge 2} E_p 
\otimes_h \max(E_q) \otimes_h E_r$. Let $C'$ be
the $cb$-norm of this map. For  clarity let us
denote again
$$ {   
\cale_{d-1}(E)  = \sum_{\scriptstyle p+q+r=d\atop
\scriptstyle q\ge 2} E_p 
\otimes_h \max(E_q)\otimes_h E_r.}$$
Then, after composing with $a_m$ and $b_m$, we obtain that the identity map from
$E^d_m = \max(\ell^m_1)_d$ to $\cale_{d-1}(\max(\ell^m_2))$ has $cb$-norm $\le 
C^d C'$.

\n But by the remark following Sublemma~2.6
(noting
$d_{cb}(\max(\ell^m_2), R_m) = 
 d_{cb}(\max(\ell^m_2), C_m) \break\le \sqrt m$),
we have
$$\|\cale_{d-1}(\max(\ell^m_2))\to
\max(\max(\ell^m_2)_d)\|_{cb} \le m^{d-2\over 
2}$$ hence a fortiori
$$\|\cale_{d-1}(\max(\ell^m_2))\to \max(\ell^m_2 \otimes_2\cdots \otimes_2 
\ell^m_2)\|_{cb} \le m^{d-2\over 2}.$$
Hence we obtain
$$\|J_{m,d}\|_{cb} \le C^d C' m^{d-2\over 2}$$
which contradicts Theorem~0.2.

\n Thus, assuming $C<\infty$, we have proved that
$d(A_d(E)) =d$.  This implies Theorem 0.3.
 Indeed, 
for any infinite dimensional  Banach space $B$,
there is, for any $\vp>0$, a factorization
$\ell^m_1  {\buildrel a_m \over {\hbox to
20pt{\rightarrowfill}}} B {\buildrel b_m \over 
{\hbox to 20pt{\rightarrowfill}}} \ell^m_2$ of
the canonical inclusion $\ell^m_1 
\to \ell^m_2$ with 
$\sup\limits_m \|a_m\|\ \|b_m\|\le 1+\vp$.
This follows from   Dvoretzky's Theorem (\cf 
e.g.\ [P4, p. 41]) applied to the dual $B^*$.
Therefore, if
$E = 
\max(B)$, we have the factorization considered in the first part of this proof 
with $C = 1+\vp$.\qed

 As already explained in the introduction, when 
$d>2$ we cannot prove Sublemma~2.8 using an independent collection of Gaussian 
random matrices indexed by $[1,\ldots, m]^d$ as a substitute for the collection 
$(U_i)$. To convince the reader of this impossibility we will now give the 
estimates resulting from this choice. Although we include them for the record, 
they may be of independent interest.

Let $Y$ be an $N\times N$ random matrix for which the entries $\{Y(i,j) \mid 
1\le i,j\le N\}$ are independent complex Gaussian random variables with 
$\EE|Y(i,j)|^2 = N^{-1}$ for all $i,j$.
To abbreviate,  we will say below
  that such a random  matrix is ``standard of size
$N\times N$''. By a well known result (which  follows from
(1.1)), there is an absolute constant $K'$ such that for all $N$
$$(\EE\|Y\|^2_{M_N})^{1/2}\le K'.\leqno (2.5)$$
More generally, by the
 concentration of measure phenomenon (see e.g.
[P4 , p. 44]), whenever 
$(Y_1,\ldots, Y_R)$ is an independent collection of copies of $Y$, and $N\ge 
\hbox{Log } R$, we have 
$$(\EE \sup_{1\le r\le R} \|Y_r\|^2)^{1/2} \le K''\leqno (2.6)$$
where $K''$ is another absolute constant. Now let $\{Y_i\mid
i\in [1,\ldots,  m]^d\}$ be an independent collection of copies
of the variable $Y$. For each $i  = (i_1,\ldots, i_d)$ in
$[1,\ldots, m]^d$, we let $e_i = e_{i_1}\otimes \cdots 
\otimes e_{i_d}$ as before and we consider the random variable
$$y = \sum_{i\in[1,\ldots, m]^d} Y_i \otimes e_i .$$ By an
abuse of notation, we will consider
$e_i$ as the natural basis either of 
$E^d_m$ or of its dual ${E^{d}_m}^*$. Thus we may
view $y$ as a random element either of 
$M_N(E^d_m)$ or of $M_N({E^{d}_m}^*)$. Recall that
if
$v\colon \ E^d_m\to M_N$ is  the linear mapping
which takes $e_i$ to $Y_i$, then
$$\|v\|_{cb(E^d_m,M_N)}  =
\|y\|_{M_N({E^{d}_m}^*)}.\leqno (2.7)$$

Then, we can state

\proclaim Theorem 2.9. There are absolute constants $K$ and $\delta>0$ such 
that: \ if $d$ is even we have
$$\leqalignno{\delta m^{d/4} &\le
(\EE\|y\|^2_{M_N({E^d_m}^*)})^{1/2} \le K  m^{d\over 4}&(2.8)\cr
\delta m^{3d/4} &\le (\EE\|y\|^2_{M_N(E^d_m)})^{1/2} \le K m^{3d\over 
4}&(2.9)}$$
while if $d\ge 3$ is odd we have
$$\leqalignno{\delta m^{d-1\over 4} &\le
(\EE\|y\|^2_{M_N(  {E^{d}_m }^*)})^{1/2} \le  Km ^{d-1\over
4}&(2.10)\cr
\delta m^{3d+1\over 4} &\le (\EE\|y\|^2_{M_N(E^d_m)})^{1/2} \le K m^{3d+1\over 
4}.&(2.11)}$$

\pf For simplicity of notation we set
$$a = (\EE\|y\|^2_{M_N( {E^{d}_m}^* )})^{1/2} \quad
\hbox{and}\quad b = (\EE 
\|y\|^2_{M_N(E^d_m)})^{1/2}.$$
We first claim that $ab\ge m^d$. Indeed, using (2.7) it is easy to check (as in 
the proof of Lemma~1.4) that
$$ab\ge \EE\left\|\sum Y_i\otimes \overline Y_i\right\| \ge \EE N^{-1}\left(tr\ 
\sum Y_iY^*_i\right) = m^d.$$
From this claim it follows that it is enough to prove both
upper bounds in (2.8)  and (2.9) (or in (2.10) and (2.11)), the
lower bounds then follow automatically.

Now to prove these upper bounds, we will use
the associativity of $\otimes_h$ and  the
(completely) isometric identities 
$$M_n(E) \simeq C_n \otimes_h E \otimes_h R_n\leqno(2.12)$$
and 
$$C_m  \otimes_h C_n= C_{mn}\quad{\rm and}\quad 
R_m  \otimes_h R_n= R_{mn}$$
proved in [BP]  and [ER2] for any operator 
space
$E$, where $C_n$ (resp.\
$R_n$) denotes the column (resp.\ row)
$n$-dimensional Hilbert  space. We have
completely contractive (in short c.c.) inclusions $C_m\to 
\ell^m_\infty$ and $R_m\to \ell^m_\infty$ (induced by the identity). Moreover 
${ E^{d}_m }^*= \ell^m_\infty \otimes_h \cdots\otimes_h
\ell^m_\infty$. Therefore,  taking  first $E$
 one dimensional in (2.12), we have  c.c.\
inclusions 
$$M_{m^{d/2}} \simeq C_{m^{d/2}} \otimes_h R_{m^{d/2}} \to
{E^{d}_m }^*\quad (d 
\hbox{ even})$$
and taking now $E=\ell^m_\infty$   in
(2.12), we have
$$M_{m^{d-1\over 2}}(\ell^m_\infty) \simeq C_{m^{d-1\over
2}} \otimes_h 
\ell^m_\infty \otimes_h R_{m^{d-1\over 2}} \to
{E^{d}_m}^*\quad (d \hbox{ odd}).$$ Hence, if $m$
is even we find
$$a \le (\EE\|y\|^2_{M_N(M_{m^{d/2}})})^{1/2}.$$
But now the matrix $m^{-d/4}y$ is ``standard'' of size $Nm^{d/2}\times 
Nm^{d/2}$ (in the above sense), hence by (2.5) $a\le K'm^{d/4}$, whence the 
right side of (2.8). Similarly, if $d$ is odd and
$d\ge 3$ we find using (2.6)  instead
$$a \le K'' m^{d-1\over 4}$$
whence the right side of (2.10).

We now turn to the upper bounds for $b$. Since the identity maps $C_m\to 
\ell^m_1$, 
$R_m\to \ell^m_1$ and $\ell^m_\infty\to \ell^m_1$ have $cb$-norms respectively 
equal to $m^{1/2}, m^{1/2}$ and $m$ we have if $m$ is even: \ $\|M_{m^{d/2}} \to 
E^d_m\|_{cb}\le m^{d/2}$ and if $m$ is odd: \ $\|M_{m^{d-1\over 
2}}(\ell^m_\infty)  \to E^d_m\|_{cb}\le m^{d+1\over 2}$. Hence repeating the 
preceding argument we find $b\le K'm^{3d/4}$ if $d$ is even and $b\le K''m^{3d+1 
\over 4}$ if $d$ is odd.
We thus obtain the right sides of (2.10) and
(2.12), which concludes the proof  of Theorem~2.9
by our original claim that $ab\ge m^d$.\qed

\n {\bf \S 3. Complements}

In this section, we wish to develop  several points which have been overlooked 
in [P1]. For the sake of generality, we return to the framework of ``similarity 
settings''. A similarity setting is a triple $(i,E,\cala)$ where $\cala$ is a 
unital algebra, $E$ is an operator space and $i\colon\ E\to \cala$ is a linear 
embedding. We will always assume (to avoid degeneracy) that there is at least 
one injective morphism $u_0\colon \ \cala\to B(H)$ with $\|u_0i\|_{cb}\le 1$. We 
will also assume that $\cala$ is generated by $i(E)$ and the unit. For any $c\ge 
1$, we denote by $\calc_c$ the class of all morphisms $u\colon\ \cala\to B(H)$ 
with $\|u i\|_{cb}\le c$ (and, say, ${\rm
card}(H)\le {\rm
card}(\cala)$) . We then define an operator space
$\tilde A_c$ as  follows:\ we introduce an
embedding
$$J_c\colon \ \cala\to \bigoplus_{u\in\calc_c} B(H_u)$$
by setting
$$J_c(x) = \bigoplus_{u\in \calc_c} u(x).\leqno \forall~x\in \cala$$
This embedding provides us with a norm on $\cala$. We denote by $\tilde A_c$ the 
completion of $\cala$ for the corresponding norm. Clearly $\tilde A_c$ is 
actually an operator algebra and (by construction) $J_c$ extends to an isometric 
morphism from $\tilde A_c$ into $B(\calh)$ with
$\calh = \bigoplus\limits_{u\in 
\calc_c} H_u$.

Let $OA(E)$ be the universal operator algebra of $E$ as defined in \S
2. Since $i\colon \ E\to \tilde A_1$ is
completely contractive it  extends to a c.c.\
morphism $\pi_1\colon \ OA(E)\to \tilde A_1$
which is a  complete metric surjection.

The next result is a reformulation of Theorem~1.7 in [P1]
(the latter was inspired by Peller's results in
[Pe]).

\proclaim Theorem 3.1. Let $c\ge 1$. Consider $f$ in $\calk\otimes \cala$ with 
$\|f\|_{\calk\otimes_{\rm min} \tilde A_c} <1$. We denote by $E^{(j)}$ the space 
$E \otimes\cdots \otimes E$ ($j$-times) viewed
 as a subspace of $\calt(E) 
\subset OA(E)$. Then for some $N\ge 1$ there are elements $F_j$ in $\calk\otimes 
E^{(j)}$ $(0 \le j \le N)$ satisfying
$$\sup_{j\ge 0} c^j \|F_j\|_{\calk \otimes_{\rm min} OA(E)} \le \left\| 
\sum\nolimits^N_0 c^jF_j\right\|_{\calk \otimes_{\rm min} OA(E)} < 1 \leqno 
(3.1)$$
such that
$$(Id_\calk \otimes \pi_1)
\left(\sum\nolimits^N_0 F_j\right) = f.$$

The next result improves on Theorem~2.5 in [P1]
(and bypasses its Lemmas~2.2 and  2.3).

\proclaim Theorem 3.2. Let $(i,E,\cala)$ be as above and let $1 \le \theta < c < 
\infty$. The following are equivalent:
\item{(i)} Every morphism in $\calc_c$ is similar to a morphism in 
$\calc_\theta$.
\item{(ii)} Every morphism $u$ in $\bigcup\limits_{b>\theta} \calc_b$ is similar 
to a morphism in $\calc_\theta$.
\item{(iii)} There are a constant $C$ and an
integer $d$ so that, for any morphism 
$u\colon \ \cala \to B(H)$ such that $ui$ is c.b.
we have $$\|u\|_{cb(\tilde A_\theta,
B(H))} \le C \sum^d_{j=0} 
\|ui\|^j_{cb}.$$
\ms

\n {\bf Proof.} By Paulsen's results (see [P1,
Prop.~1.8] for details), it is the same to say
that
$u$ is similar to a morphism in $\calc_\theta$,
or that $\|u\|_{cb(\tilde A_\theta,
B(H))}<\infty$. Thus it is clear that (iii)
implies (ii), and trivial that (ii)
implies (i). Thus, it suffices to show that
(i)
implies (iii).
Again by Paulsen's results (see [P1,
Prop.~1.8] for details), 
  (i) holds 
iff the canonical morphism $\tilde A_c\to \tilde A_\theta$ is a complete 
isomorphism, \ie there is a constant $K>0$ such that for any $f$ in 
$\calk\otimes \cala$ we have
$$\|f\|_{\calk\otimes_{\rm min} \tilde A_c} \le K\|f\|_{\calk \otimes_{\rm min} 
\tilde A_\theta}.\leqno(3.2)$$
Assume that this holds.  Then  select the smallest
integer $d$ such that $\sum\limits_{j>d}
\left({\theta \over  c}\right)^j \le 1/2K$. We
will show that (iii) follows for some $C$. 
Let $T_d$ be the closed subspace of $OA(E)$ generated by  $\CC \oplus E \oplus 
E^{\otimes 2} \oplus\cdots \oplus E^{\otimes d}$. Let $u\colon \ \cala\to B(H)$ 
be a morphism and let $b = \|ui\|_{cb}$. 

\n Note that by [P1, Prop.~1.10] and by
Proposition~2.1 (i) above, we have
$$\|u\pi_{1\mid T_d}\|_{cb} \le \sum^d_{j=0} b^j.\leqno (3.3)$$
Now consider $f$ in $\calk \otimes \cala$ with $\|f\|_{\calk\otimes_{\rm min} 
\tilde A_\theta} < 1$ and hence by (3.2) $\|f\|_{\calk\otimes_{\rm min} \tilde 
A_c} < K$.

\n We claim that $f$ can be decomposed in $\calk
\otimes\cala$ as $f = (I_\calk 
\otimes \pi_1) (x_0) + f'$ with $x_0\in \calk \otimes T_d$, $f'\in \calk \otimes 
\cala$ satisfying
$$\|x_0\|_{\calk \otimes_{\rm min}T_d} \le C'\quad \hbox{and}\quad \|f'\|_{\calk 
\otimes_{\rm min}\tilde A_\theta} < {1\over 2}.
$$ where $C' = K\sum^d_0 c^{-j}$.

\n From this claim (iii) follows immediately.
Indeed, iterating the claim, we find  a sequence
$(x_n)_{n\ge 0}$ in $\calk \otimes  T_d$ such
that $\|x_n\|_{\calk 
\otimes_{\rm min} T_d} \le C' 2^{-n}$ and
$$\left\|f - (I_\calk \otimes \pi_1)
\left(\sum\nolimits^n_0 
x_j\right)\right\|_{\calk \otimes_{\rm min}\tilde
A_\theta} < 2^{-n-1}.$$ Let $y$ be the sum of the
series $\sum^\infty_0 x_j$ which converges in
$\calk 
\otimes_{\min} T_d$. Note that $\|y\|_{\calk
\otimes_{\rm min} T_d} \le 2C'$.
Since $\pi_{1\mid T_d} \colon \ T_d\to \tilde
A_1$ actually is continuous into 
$\tilde A_\theta$, we can write $f = (I_\calk \otimes \pi_1)(y)$ in $\calk 
\otimes_{\rm min} \tilde A_\theta$. Therefore
$$\leqalignno{\|(I_\calk\otimes u)(f)\| &\le \|u\pi_{1\mid T_d}\|_{cb} 
\|y\|_{\calk \otimes_{\rm min} T_d}\cr
&\le \left(\sum\nolimits^d_0 b^j\right) 2C',&\hbox{hence by (3.3)}}$$
which implies that (iii) holds with $C=2C'$. Thus, to complete the proof, it 
suffices to prove the above claimed decomposition for $f$. By Theorem~3.1, we 
can write
$$f = \sum\nolimits^N_0 (I_\calk \otimes \pi_1)(F_j)$$
with $ \sup_{j\ge 0} {c^j\left\| 
F_j\right\|_{\calk
\otimes_{\rm min} OA(E)} }< K$. We can  always
assume (adding zero terms if necessary) that
$N\ge d$. Then we have
$$\leqalignno{f &= (I_\calk \otimes \pi_1)(x_0) + f'\cr
x_0 &= \sum\nolimits^d_0 F_j&\hbox{where}\cr
f' &= \sum_{j>d} (I_\calk \otimes \pi_1)(F_j).&\hbox{and}}$$
Then by (3.1) and our original choice of $d$, we
have
$$\eqalign{\|f'\|_{\calk\otimes_{\rm min} \tilde
A_\theta} &\le \sum_{j>d} 
\theta^j \|F_j\|_{\calk \otimes_{\rm min} OA(E)}\cr
&< \sum_{j>d} (\theta/c)^j K \le 1/2.}$$
On the other hand, by (3.1) again
$$\|x_0\|_{\calk \otimes_{\rm min} T_d} \le \sum\nolimits^d_0 \|F_j\|_{\calk 
\otimes_{\rm min} OA(E)} \le C',$$
which establishes the above claim.\qed

\n {\bf Remark 3.4.} Just like in Theorem~2.6 in [P1] the preceding proof works 
just as well if we replace $\calk$ throughout the proof by a subspace $X\subset 
\calk$ for which there is a projection $P\colon \ \calk \to X$ with 
$\|P\|_{cb}=1$. If $X \otimes_{\rm min} \tilde A_\theta$ is isomorphic to $X 
\otimes_{\rm min} \tilde A_c$, then $X \otimes_{\rm min} \tilde A_\theta = X 
\otimes_{\rm min} \tilde A_b$ for all $b\ge \theta$. In particular this applies 
when $X$ is 1-dimensional. In this case, we find that if $\tilde A_\theta$ and 
$\tilde A_c$ have equivalent norms, then $\tilde A_\theta$ and $\tilde A_b$ have 
equivalent norms for all $b\ge \theta$.\ms

\n {\bf Proof of Theorems 0.4 and 0.5.} These
statements are nothing but  Theorem~3.2 in the
particular case $E = \max(A)$ with $i$ equal to
the identity  on $A$. \qed \ms

\n {\bf Remark.} We refer the interested reader to
[P7] for more information on  the themes of the
present paper in the context of uniformly bounded
group  representations on locally compact groups.
(The presentation of [P7] assumes  very little
familiarity with operator spaces.)\ms

\n {\bf Remark.} Theorem 3.1 can be applied in the situation considered in [P2]. 
Let $A$ be a unital operator algebra and let $A_1,A_2$ be unital (closed) 
subalgebras, let $\cala$ be the algebra generated by $A_1\cup A_2$. We assume 
$\cala$ dense in $A$. The associated similarity setting is: \ $E = A_1\oplus_1 
A_2$ (operator space $\ell_1$-direct sum) with $i\colon \ E\to \cala$ defined by 
$i((x_1,x_2)) = x_1+x_2$.

\n Let us denote here for simplicity
$${\calk}(A) = \calk \otimes_{\rm min} A.$$
Clearly ${\calk}(A)$ is an operator algebra which we may view as formed of bi-infinite 
matrices with entries in $A$.

We will say that $(A_1,A_2)$ generate $A$ with length $\le d$ if any $x$ in 
${\calk}(A)$ can be written as
$$x = x_1x_2 \ldots x_d + y_1y_2\ldots y_d\leqno
(3.4)$$ with $x_i\in {\calk}(A_1)\cup {\calk}(A_2)$, $y_i \in
{\calk}(A_1) \cup {\calk}(A_2)$ $(1\le i \le d)$  and also
$x_1\in {\calk}(A_1)$ and $y_1\in {\calk}(A_2)$. This implies
that the natural  product map is a surjection
from $[{\calk}(A_1) \widehat\otimes \cdots ] \oplus_1
\break [{\calk}(A_2) \widehat\otimes\cdots]$ onto
${\calk}(A)$. Hence by the open mapping theorem  (and
by a well known ``matrix trick''), there is a
constant $K$ such that we can  always find
$x_i,y_i$ as above satisfying moreover
$$\prod^d_{i=1} \|x_i\| + \prod^d_{i=1} \|y_i\| \le K\|x\|.$$
We denote by $\ell(A_1,A_2)$ the smallest $d$ such that $\ell(A_1,A_2)\le d$.

Note that this definition is equivalent to [P2, Definition~5]:\ indeed an 
elementary argument allows to pass from the approximate version of
(3.4) given  in [P2] to equality in (3.4). In
this case, $\ell(A_1,A_2)$ is equal to the 
degree of the setting $(i,E,\cala)$, and there is
a one to one correspondence  between  morphisms
$u\colon \ \cala \to B(H)$ with $\|ui\|_{cb}\le
c$ and pairs of  morphisms $u_i\colon \ A_i\to
B(H)$ with $\max\limits_{i=1,2} \{\|u_i\|_{cb}\} 
\le c$. We refer the reader to [P2] for more variations on this theme.

\vfill\eject

\centerline {\bf References}

\item{[BL]} J. Bergh and J. L\"ofstr\"om. Interpolation spaces. An
 introduction. Springer Verlag, New York. 1976.

%\item{[BRS]} D. Blecher, Z. J. Ruan  and  A.
%Sinclair.  A characterization of operator
%algebras. J. Funct. Anal. 89 (1990) 188-201.

\item{[B]} D. Blecher.  Tensor products of
operator spaces II. Can. J. math.  44 (1992)
75-90.

\item{[BP]} D. Blecher and V. Paulsen.  Tensor products of
operator spaces   J. Funct. Anal.  99 (1991) 262-292.

%\item{[BP2]}  $\underline{\hskip1.5in}$. 
%Explicit construction of universal operator
%algebras and applications to polynomial
%factorization. Proc. Amer. Math. Soc.   112
%(1991) 839-850. 

\item{[BRS]} D. Blecher, Z. J. Ruan and A. Sinclair.
 A characterization
of operator
algebras. J. Funct. Anal. 89 (1990) 188-201.

\item{[Bo]} M. Bo\.{z}ejko. 	
Positive-definite kernels, length functions on
groups and a noncommutative von Neumann
inequality.   Studia Math.   95 (1989) 107-118.

\item{[CS1]} E. Christensen and A. Sinclair. 
Representations of completely bounded multilinear operators.
J. Funct. Anal. 72 (1987) 151-181.

\item{[CS2]}   $\underline{\hskip1.5in}$.  A
survey of completely bounded operators.  Bull.
London Math. Soc.  21 (1989) 417-448.
 
 \item{[ER1]} E. Effros     and Z.J. Ruan.  A new
approach to operators spaces.
 Canadian Math. Bull. 
34 (1991) 329-337.

\item{[ER2]} $\underline{\hskip1.5in}$. Self
duality
 for the Haagerup
 tensor product and Hilbert space factorization.  J. Funct. Anal. 100
(1991) 257-284.

 \item{[ER3]} $\underline{\hskip1.5in}$. Mapping
spaces and liftings for operator spaces.  Proc.
London Math. Soc.  69 (1994) 171-197.

\item{[ER4]} $\underline{\hskip1.5in}$. On
approximation properties for operator spaces,
International J. Math. 1 (1990) 163-187.

\item{[Ge]} S. Geman. A limit theorem for the
norm of random matrices. Ann. Prob. 8 (1980)
252-261.

 \item{[HT]} U. Haagerup and S. Thorbj{\o}rnsen.
Random matrices and $K$-theory for exact
$C^*$-algebras. Odense University preprint 1998.
 
\item{[Ka]} R.  Kadison  On the orthogonalization
of operator representations. Amer. J. Math.  77
(1955) 600-620.

 \item{[MP]} M. Marcus and G. Pisier. {  Random
Fourier series with Applications to Harmonic
Analysis.}  Annals
of Math. Studies n$^\circ$101, Princeton
University Press. (1981). 

\item{[Pa1]} V. Paulsen.   Completely bounded maps and
dilations.  Pitman Research Notes in
Math. 146, Longman, Wiley, New York, 1986.

\item{[Pa2]}   $\underline{\hskip1.5in}$. The
maximal operator space of a normed space.
 Proc. Edinburgh Math. Soc.  39 (1996) 309-323.

\item{[Pa3]}   $\underline{\hskip1.5in}$.
 Representation of
Function algebras, Abstract
 operator spaces and Banach space Geometry. J.
Funct. Anal. 109 (1992) 113-129.

 \item{[PS]} V. Paulsen and R. Smith. Multilinear
maps and tensor norms on operator systems. J.
Funct. Anal. 73 (1987) 258-276.

 \item{[Pe]} V. Peller.  Estimates of functions
 of power bounded operators on Hilbert space. J.
 Oper. Theory  7 (1982) 341-372.

 \item{[Pes]} V. Pestov.  Operator spaces and residually
finite-dimensional $C^*$-algebras.
J. Funct. Anal. 123 (1994) 308-317.  

\item{[P1]}  G. Pisier.  The similarity degree of
an operator algebra.
   St. Petersburg Math. J. 10 (1999) 103-146.
 
\item{[P2]}   $\underline{\hskip1.5in}$. Joint
similarity problems and the generation of
operator algebras with bounded length. Integr.
Equ. Op. Th.
 31 (1998) 353-370.

\item{[P3]}   $\underline{\hskip1.5in}$.
Similarity problems and completely bounded maps.
Springer Lecture notes 1618 (1995).
 
\item{[P4]} $\underline{\hskip1.5in}$. The volume of Convex Bodies and Banach
Space Geometry.  Cambridge University
Press, 1989.

\item{[P5]} $\underline{\hskip1.5in}$. The
operator Hilbert space $OH$, complex
interpolation and tensor norms. Memoirs Amer.
Math. Soc.  vol. 122 , 585 (1996) 1-103.

\item{[P6]} $\underline{\hskip1.5in}$. Noncommutative
 vector valued $L_p$-spaces and completely
$p$-summing maps. Ast\'erisque (Soc. Math. France) 247 (1998) 1-131.

\item{[P7]} $\underline{\hskip1.5in}$. Are
unitarizable groups amenable ? preprint, 1999.

%\item{[V]}  Varopoulos N.  On an inequality of
%von Neumann and an application of the metric
%theory of tensor products to Operators Theory. 
%J. Funct. Anal.  16 (1974), 83-100.

 \bye
\vskip12pt

Texas A\&M  University

College Station, TX 77843, U. S. A.

and

Universit\'e Paris VI

Equipe d'Analyse, Case 186,
 
75252 Paris Cedex 05, France
\end